\documentclass[11pt]{article}
\usepackage{ytableau}
\usepackage{bm}
\usepackage{amssymb,array}
\usepackage[margin=1in]{geometry}
\usepackage{mathrsfs}
\usepackage{stmaryrd}
\usepackage{amsfonts,amsmath,amssymb,amscd,amsthm}
\usepackage{shadow}
\usepackage{epsfig,epstopdf}
\usepackage{graphicx}
\usepackage{color}
\usepackage{multirow}
\usepackage{booktabs} 
\usepackage{enumerate}
\usepackage{longtable}
\usepackage{pstricks,multido}
\usepackage[colorlinks,urlcolor=purple,linkcolor=red,anchorcolor=green,citecolor=blue]{hyperref}
\usepackage{qtree}
\usepackage{pxfonts}
\allowdisplaybreaks
\newtheorem{thm}{Theorem}[section]
\newtheorem{lem}[thm]{Lemma}

\newtheorem{remark}{Remark}[section]
\newtheorem{example}{Example}[section]

\newtheorem{conj}{Conjecture}[section]

\makeatletter
\newcommand{\neutralize}[1]{\expandafter\let\csname c@#1\endcsname\count@}
\makeatother

\newenvironment{conjprime}[1]
  {\renewcommand{\theconj}{\ref*{#1}$'$}%
   \neutralize{conj}\phantomsection
   \begin{conj}}
  {\end{conj}}

\newcommand{\tc}[1]{\multicolumn{1}{c}{#1}} 

\numberwithin{figure}{section}
\numberwithin{table}{section}

\def\qed{\hfill \rule{4pt}{7pt}}

\def\pf{\noindent {\it{Proof.} \hskip 2pt}}

\numberwithin{equation}{section}
\begin{document}
\begin{center}
{\Large\bf  Combinatorial proofs and generalizations  of conjectures \\[3pt]
  related to Euler's partition theorem  }
\end{center}

\begin{center}

{ Jane Y.X. Yang }

{\small
School of Science\\
Chongqing University of Posts and Telecommunications\\
Chongqing 400065, P.R. China

yangyingxue@cqupt.edu.cn}

\end{center}

\begin{abstract}
In a recent work,  Andrews
gave analytic proofs of two conjectures concerning some variations of two combinatorial identities
between partitions of a positive integer into odd parts and partitions into distinct parts discovered
by Beck. Subsequently, using the same method as Andrews,  Chern  presented the analytic proof of another Beck's conjecture
relating the gap-free partitions and distinct partitions with  odd length.
However, the combinatorial interpretations of these conjectures
are still unclear and required.
In this paper, motivated by Glaisher's bijection,
we give the combinatorial proofs of these three conjectures
directly or by proving more generalized results.
\end{abstract}



\section{Introduction}

A \emph{partition} \cite{Andrews98} of $n$ is a finite nonincreasing sequence of positive integers
$(\lambda_1, \lambda_2, \ldots, \lambda_\ell)$ such that $n=\lambda_1+\lambda_2+\cdots+\lambda_\ell$.
We write
$\lambda=(\lambda_1,\lambda_2,\ldots,\lambda_\ell)$
and call $\lambda_i$'s the \emph{parts} of $\lambda$.
If a part $i$ has multiplicity  $m_i$ for $i\geq 1$, we also write $\lambda$ as $(1^{m_1},2^{m_2},\ldots)$,
where the superscript $m_i$ is neglected provided $m_i=1$.
The  \emph{weight} of $\lambda$ is the sum of all parts, which is denoted by $|\lambda|$,
and the \emph{length }of $\lambda$ is the number of parts,
which is denoted by $\ell(\lambda)$.
The \emph{conjugate} of $\lambda$ is the partition
$\lambda'=(\lambda'_1,\lambda'_2,\ldots,\lambda'_{\lambda_{1}})$,
where
$\lambda'_{i}=| \{\lambda_{j}:\lambda_{j}\geq i,1\leq j \leq \ell \}|$ for $1\leq i \leq \lambda_{1}$,
or $\lambda'$ can be equivalently expressed as $(1^{\lambda_1-\lambda_2},2^{\lambda_2-\lambda_3},\ldots,\ell-1^{\lambda_{\ell-1}-\lambda_{\ell}},\ell^{\lambda_{\ell}})$.
A partition $\lambda=(\lambda_1,\lambda_2,\ldots,\lambda_\ell)$
is called a distinct partition if $\lambda_1>\lambda_2>\cdots>\lambda_\ell$,
and an odd partition if $\lambda_i$ is odd for all $1\leq i\leq \ell$, respectively.
In 1748,
by using generating functions,
Euler \cite{Euler} gave the celebrated partition theorem as follows.
\begin{thm}[Euler's partition theorem]
The number of distinct partitions of $n$ equals the number of odd partitions of $n$.
\end{thm}

After Euler's partition theorem was proposed, there have been many extensions and
refinements,
the famous ones of which are
Glaisher's theorem and Franklin's theorem,
and the reader can refer to \cite{Pak,Remmel,Wilf} for more details.
Glaisher \cite{Glaisher}  bijectively proved  the following extension.
\begin{thm}[Glaisher's theorem]\label{glaisherthm}
   For any positive integer $k$, the number of partitions of $n$ with no part occurring $k$ or more times
  equals the number of partitions of $n$
  with no part divisible by $k$.
\end{thm}
In 1882, Franklin \cite{Franklin,Sylvester} acquired a more generalized result by giving constructive proof of the following theorem. 
Franklin \cite[p. 268]{Sylvester} also asserted that the generating functions are easily obtained.
\begin{thm}[Franklin's theorem]\label{franklinthm}
   For any positive integer $k$ and nonnegative integer $m$,
  the number  of partitions of $n$
  with $m$ distinct parts
  each occurring $k$ or more times
  equals
  the number of partitions of $n$ with exactly $m$ distinct parts divisible by $k$.
\end{thm}
Thus by taking $m=0$, Franklin's theorem degenerates to Glaisher's theorem,
then by taking $k=2$, Glaisher's theorem
gives Euler's partition theorem.

From the works of Andrews
and Chern,
we noticed three conjectures  posed by Beck
concerning  some variations of odd partitions
and distinct partitions,
which were only analytically proved by Andrews \cite{Andrews} and Chern \cite{Chern}
via differentiation technique in $q$-series introduced by Andrews \cite{Andrews}.
In this paper, by extending Glaisher's bijection,
we give the combinatorial proofs of these three conjectures
directly or by proving more generalized results.
For consistency of notations,
we utilize the same notations in \cite{Andrews} and \cite{Chern}
in the rest of paper as far as possible.

Let $a(n)$ denote the number of partitions
of $n$ with only one even part which  is possibly repeated.
Beck \cite{Beck1} proposed the following conjecture:
\begin{conj}\label{conj1}
$a(n)$ is also the difference between the number of parts in the  odd partitions of $n$ and the number of parts in the distinct partitions of $n$.

\end{conj}


Let $c(n)$ denote the number of partitions of $n$ in which exactly one part is repeated.
Let $b(n)$ be the difference between the number of parts in the odd partitions of $n$ and the number of parts in the distinct partitions of $n$.
Andrews \cite{Andrews} analytically proved the following theorem by differentiation technique in $q$-series,
which confirms the conjecture posed by Beck \cite{Beck1}:
\begin{thm}[{\cite[Theorem 1]{Andrews}}]\label{thmconj1andrews}
For all $n\geq 1, a(n)=b(n)=c(n)$.
\end{thm}
Later Fu and Tang \cite[Theorem 1.5]{Fu} extended
the result of Andrews and gave the analytic proof.

For $k\geq1$, let $\mathcal{O}_k(n)$ be
the set of partitions of $n$
with no part divisible by $k$
and $\mathcal{D}_k(n)$ be the set of
partitions of $n$ with no  part occurring  $k$ or more times, respectively.
Let $\mathcal{O}_{1,k}(n)$ be the partitions of $n$
with exactly one part (possibly repeated) divisible by $k$.
As one of main results in this paper, the following theorem generalizes Conjecture \ref{conj1}.
\begin{thm}\label{thmofconj1}
For $k\geq2$ and $n\geq0$, we have
\begin{equation*}
  \left|\mathcal{O}_{1,k}(n)\right|=\frac{1}{k-1}\cdot\Bigg(\sum_{\lambda\in\mathcal{O}_k(n)}\ell(\lambda)-\sum_{\lambda\in\mathcal{D}_k(n)}\ell(\lambda)\Bigg).
\end{equation*}
\end{thm}
Thus, letting $k=2$ in Theorem \ref{thmofconj1} reduces the set $\mathcal{O}_{1,k}(n)$
to the set of  partitions of $n$ with only one even part which is possibly repeated,
and the set $\mathcal{O}_k(n)$ (resp. $\mathcal{D}_k(n)$) to the set of odd (resp. distinct) partitions of $n$,
which  gives the combinatorial proof of Conjecture \ref{conj1}.

Let $a_1(n)$ denote the number of partitions of $n$ such that there is exactly one part occurring three times while all other parts occur only once.
Beck \cite{Beck1} made the following conjecture:
\begin{conj}\label{conj2}
$a_1(n)$ is also the difference between the number of parts in the distinct partitions of $n$ and the number of distinct parts in the odd partitions of $n$.
\end{conj}


Let $b_1(n)$ be the  difference between the total number of parts
in the partitions of $n$ into distinct parts
and the total number of distinct parts
in the partitions of $n$ into odd parts.
This conjecture was also proved by Andrews \cite{Andrews} with analytic method.

\begin{thm}[{\cite[Theorem 2]{Andrews}}]\label{andrewsthm2}
  $a_1(n)=b_1(n).$
\end{thm}

Denote by $\mathcal{T}_k(n)$ the set of partitions of $n$ such that
there is exactly one part occurring more than $k$ times and less than $2k$ times
while all other parts occur less than $k$ times.
Motivated by Glaisher's bijection,
instead of  proving Conjecture \ref{conj2} directly,
we prove a more generalized theorem as follows.
\begin{thm}\label{thmofconj2}
  For $k\geq1$ and $n\geq 0$, we have
  \begin{equation*}\label{eqofthmofconj2}
    \left|\mathcal{T}_k(n)\right|=\sum_{\lambda\in\mathcal{D}_k(n)}\bar{\ell}(\lambda)-\sum_{\lambda\in\mathcal{O}_k(n)}\bar{\ell}(\lambda),
  \end{equation*}
  where  $\bar{\ell}(\lambda)$ is the number of distinct parts in $\lambda$.
\end{thm}
Thus by taking $k=2$, $\mathcal{T}_2(n)$ becomes the set of partitions of $n$
with exactly one part occurring three times and other parts occurring only once while $\mathcal{D}_2(n)$ and $\mathcal{O}_2(n)$ become the set of distinct and odd partitions of $n$, respectively, which gives a positive answer to Conjecture \ref{conj2}.


A partition $\lambda=(\lambda_1,\ldots,\lambda_\ell)$ is called a \emph{gap-free}, or \emph{compact}, partition if $0\leq\lambda_i-\lambda_{i+1}\leq 1$ for all $1\leq i\leq \ell-1$.
Let $a_2(n)$ be the number of gap-free partitions of $n$.
Andrews \cite{Andrews17} gave the generating function of $a_2(n)$.
Beck \cite{Beck2} proposed the following conjecture of $a_2(n)$:

\begin{conj}\label{conj3}
$a_2(n)$ is also the sum of the smallest parts in the distinct partitions of $n$ with an odd number of parts.
\end{conj}


Let $b_2(n)$ denote the sum of the smallest parts in the distinct partitions of $n$ with an odd number of parts.
Chern \cite{Chern} proved this conjecture analytically by $q$-series based on method used by Andrews in \cite{Andrews}.

\begin{thm}[{\cite[Theorem 1.2]{Chern}}]
  For all $n\geq1$, $a_2(n)=b_2(n)$.
\end{thm}

In this paper, we not only give the
combinatorial proof of Conjecture \ref{conj3},
but also study the relationship
between the gap-free partitions
and the distinct partitions with even length,
which leads us to rediscover a classical
combinatorial identity found by Fokkink--Fokkink--Wang \cite{Fokkink-Fokkink-Wang} combinatorially and
proved by Andrews \cite{Andrews08} analytically.

The rest of this paper is organized as follows.
As our main tool to prove Theorem
\ref{thmofconj1} and Theorem \ref{thmofconj2},
we give a detailed introduction of
Glaisher's bijection in Section \ref{secofGlaisher}.
In Section \ref{secconj1} and Section \ref{secconj2}, we give the combinatorial proofs of Theorem \ref{thmofconj1}
and Theorem \ref{thmofconj2}, respectively.
The combinatorial proof of Conjecture \ref{conj3} (Theorem \ref{thmofconj3}) and a similar theorem (Theorem \ref{thmofconj3'})
connecting the gap-free partitions with
the distinct partitions of even length
are contained in Section \ref{secconj3}.

\section{Glaisher's bijection}\label{secofGlaisher}

In this section, we mainly recall Glaisher's bijection \cite{Glaisher} that gives the proof of $|\mathcal{O}_k(n)|=|\mathcal{D}_k(n)|$,
since it will be used frequently throughout this paper.

Glaisher's bijection $\phi$
is defined as follows,
and one can refer to \cite{Glaisher, Pak, Remmel, Wilf} for more details.
Let $\lambda=(\lambda_1,\lambda_2,\ldots,\lambda_\ell)$
be a partition in $\mathcal{D}_k(n)$,
then for $1\leq i\leq\ell$, each part $\lambda_i$ can be uniquely written as $k^{m_i}f_i$ with $k \nmid f_i$.
Thus $\phi(\lambda)$ is
established from $\lambda$ by
replacing the part $\lambda_i$ by $k^{m_i}$
parts $f_i$ for $1\leq i\leq \ell$.
Since each $f_i$ is not divisible by $k$,
then we have $\phi(\lambda)\in\mathcal{O}_k(n)$.

\begin{example}
  If $k=3$ and $\lambda=(2^2,6,8^2,9,12)\in\mathcal{D}_3(47)$,
  then $2=3^0\cdot2$, $6=3^1\cdot2$, $8=3^0\cdot8$, $9=3^2\cdot1$ and $12=3^1\cdot4$.
  Hence $\phi(\lambda)=(1^9,2^5,4^3,8^2)\in\mathcal{O}_3(47)$.
\end{example}

In another direction, given a partition $\mu=(\mu_1,\mu_2,\ldots,\mu_\ell)\in\mathcal{O}_k(n)$,
assume that there are $m_i$ parts $i$ in $\mu$ for $1\leq i\leq n$.
Then for each $m_i\geq1$,
$m_i$ can be uniquely expressed as
$$m_i=b_{i_1}k^{a_{i_1}}+b_{i_2}k^{a_{i_2}}+\cdots+b_{i_{p_i}}k^{a_{i_{p_i}}},$$
where $a_{i_1}>a_{i_2}>\cdots>a_{i_{p_i}}\geq0$ and $1\leq b_{i_j} \leq k-1$ for $1\leq j\leq p_i$.
Hence we can construct $\varphi(\mu)$ from $\mu$ by replacing   $m_i$ parts $i$ by $b_{i_1}$ parts $i\cdot k^{a_{i_1}}$,
$b_{i_2}$ parts $i\cdot k^{a_{i_2}}$, $\ldots$ , $
b_{i_{p_i}}$ parts $i\cdot k^{a_{i_{p_i}}}$
whenever $m_i\geq1$.
It is clear that $\varphi(\mu)\in\mathcal{O}_k(n)$
because $1\leq b_{i_j}\leq k-1$ for $1\leq i\leq n$  and  for $1\leq j\leq p_i$,
and $i\cdot k^{r}=j\cdot k^s$
if and only if both $i=j$ and $r=s$ hold
since $i$ and $j$ are not divisible by $k$.
Actually, we can describe the process of
producing $\varphi(\mu)$  in a simpler way.
It is easy to check that $\varphi(\mu)$
is obtained from $\mu$ by merging  $k$ equal parts into one part
until there is no part occurring $k$ or more times.

\begin{example}
  If $k=3$ and $\mu=(1^9,2^5,4^3,8^2)\in\mathcal{O}_3(47)$,
  then $9=1\cdot3^2$, $5=1\cdot3^1+2\cdot3^0$, $3=1\cdot3^1$ and $2=2\cdot3^0$.
  Thus
  $\varphi(\mu)=(2^2,6,8^2,9,12)\in\mathcal{D}_3(47)$.
\end{example}

Therefore, the maps  $\phi\colon\,\mathcal{D}_k(n)\rightarrow\mathcal{O}_k(n)$
and
$\varphi\colon\,\mathcal{O}_k(n)\rightarrow\mathcal{D}_k(n)$
are well-defined bijections inverse to each other: $\varphi=\phi^{-1}$,
which gives Glaisher's bijection.
In the rest of paper,
according to the specific circumstances,
we are free to choose Glaisher's bijection as
$\phi$ or $\varphi$.

\section{Combinatorial proof of Theorem \ref{thmofconj1}}\label{secconj1}


In this section, we use Glaisher's bijection $\varphi\colon\,\mathcal{O}_k(n)\rightarrow \mathcal{D}_k(n)$.
Let $m$ be a positive integer,
then denote by $p(m)$  the sum of nonzero digits
in the $k$-adic representation of $m$ and
{$a(m)$ the highest exponent in the $k$-adic representation of $m$.
Precisely, $m$ can be uniquely written as $b_1k^{a_1}+b_2k^{a_2}+\cdots+b_pk^{a_p}$
with $a_1>a_2>\cdots>a_p\geq0$ and $1\leq b_i \leq k-1$ for $1\leq i\leq p$,
then $p(m)=\sum_{i=1}^pb_i$ and $a(m)=a_1$.

\begin{lem}\label{lemofm-pm}
  Let $\lambda=(1^{m_1},2^{m_2},\ldots)\in\mathcal{O}_k(n)$,
  then we have

  \begin{equation}\label{eqofm-pm}
    \ell(\lambda)-\ell(\varphi(\lambda))=\sum_{i}(m_i-p(m_i)).
  \end{equation}
\end{lem}
\pf
Assume that $\lambda$ contains $m_i$ parts $i$,
then by the construction of Glaisher's bijection $\varphi\colon\,\mathcal{O}_k(n)\rightarrow\mathcal{D}_k(n)$,
for each $m_i$,
we rewrite $m_i$ as  $b_{i_1}k^{a_{i_1}}+b_{i_2}k^{a_{i_2}}+\cdots+b_{i_{p_i}}k^{a_{i_{p_i}}}$ 
and replace $m_i$ parts $i$  by $b_{i_1}$ parts $i\cdot k^{a_{i_1}}$,
$b_{i_2}$ parts $i\cdot k^{a_{i_2}}$, $\ldots$,
$b_{i_{p_i}}$ parts $i\cdot k^{a_{i_{p_i}}}$.
Therefore the number of parts of $\varphi(\lambda)$ is decreased to $\sum_{i}p(m_i)$,
which implies \eqref{eqofm-pm} since the number of parts of $\lambda$ is $\sum_{i}m_i$.
\qed

Next we will construct a series of subsets of $\mathcal{O}_{\lambda,k,i}\subseteq\mathcal{O}_{1,k}(n)$ for any $\lambda\in\mathcal{O}_k(n)$.
Let $\lambda=(1^{m_1},2^{m_2},\ldots)\in\mathcal{O}_k(n)$.
Then for each $m_i$,
we establish a set of  partitions $\pi^i_{j,r}$ from $\lambda$ by
replacing $r\cdot k^j$ parts $i$ by
$r$ parts $i\cdot k^j$,
where $1\leq j\leq a(m_i)$
and $1\leq r\leq \lfloor m_i/k^j \rfloor$.
Thus there exists only one part $i\cdot k^j$ divisible by $k$
occurring $r$ times in $\pi^i_{j,r}$,
which implies $\pi^i_{j,r}\in\mathcal{O}_{1,k}(n)$.
Define
$$\mathcal{O}_{\lambda,k,i}=\left\{\pi^{i}_{j,r}\colon\,1\leq r\leq \lfloor m_i/k^j \rfloor \mbox{ and } 1\leq j\leq a(m_i)\right\}.$$
 Note that $\mathcal{O}_{\lambda,k,i}=\emptyset$ when $m_i<k$ since $a(m_i)=0$.
%

\begin{lem}\label{lemmaofElambdai}
Let $\lambda=(1^{m_1},2^{m_2},\ldots)\in\mathcal{O}_k(n)$, then we have
  $$\left|\mathcal{O}_{\lambda,k,i}\right|=\frac{m_i-p(m_i)}{k-1}.$$
\end{lem}
\pf
Since $m_i<k$ is trivial, we assume  $m_i\geq k$ and $$m_i=b_{i_1}k^{a_{i_1}}+b_{i_2}k^{a_{i_2}}+\cdots+b_{i_{p_i}}k^{a_{i_{p_i}}},$$
where $a_{i_1}>a_{i_2}>\cdots>a_{i_{p_i}}\geq0$,
then $p(m_i)=\sum_{j=1}^{p_i}b_{i_j}$ and $a(m_i)=a_{i_1}$.
Thus we see that
\begin{align*}
 \left|\mathcal{O}_{\lambda,k,i}\right|
       &=\left\lfloor\frac{m_i}{k}\right\rfloor+\left\lfloor\frac{m_i}{k^2}\right\rfloor+\cdots+\left\lfloor\frac{m_i}{k^{a_{i_1}}}\right\rfloor\\[6pt]
       &=\left\lfloor b_{i_1}k^{a_{i_1}-1}+b_{i_2}k^{a_{i_2}-1}+\cdots+b_{i_{p_i}}k^{a_{i_{p_i}}-1}\right\rfloor+\left\lfloor b_{i_1}k^{a_{i_1}-2}+b_{i_2}k^{a_{i_2}-2}+\cdots+b_{i_p}k^{a_{i_{p_i}}-2}\right\rfloor+\\[6pt]
       &\qquad\cdots+\left\lfloor b_{i_1}+b_{i_2}k^{a_{i_2}-a_{i_1}}+\cdots+b_{i_{p_i}}k^{a_{i_{p_i}}-a_{i_1}}\right\rfloor\\[6pt]
       &=\left(b_{i_1} k^{a_{i_1}-1}+b_{i_1} k^{a_{i_1}-2}+\cdots+b_{i_1}\right)+\left(b_{i_2} k^{a_{i_2}-1}+b_{i_2} k^{a_{i_2}-2}+\cdots+b_{i_2}\right)+\\[6pt]
       &\qquad\cdots+\left(b_{i_{p_i}} k^{a_{i_{p_i}}-1}+b_{i_{p_i}} k^{a_{i_{p_i}}-2}+\cdots+b_{i_{p_i}}\right)\\[6pt]
       &=b_{i_1}\cdot\frac{k^{a_{i_1}}-1}{k-1}+b_{i_2}\cdot\frac{k^{a_{i_2}}-1}{k-1}+\cdots+b_{i_{p_i}}\cdot\frac{k^{a_{i_{p_i}}}-1}{k-1}\\[6pt]
       &=\frac{1}{k-1}\cdot\left(b_{i_1} k^{a_{i_1}}+b_{i_2} k^{a_{i_2}}+\cdots+b_{i_{p_i}} k^{a_{i_{p_i}}}-(b_{i_1}+b_{i_2}+\cdots+b_{i_{p_i}})\right)\\[6pt]
       &=\frac{m_i-p(m_i)}{k-1},
\end{align*}
which completes the proof.
\qed

\begin{thm}\label{thmofconj1'}
   For $k\geq2$ and $n\geq0$, we have
 $$\mathcal{O}_{1,k}(n)=\bigcup_{\lambda\in\mathcal{O}_k(n)}\Bigg(\bigcup_{i}\mathcal{O}_{\lambda,k,i}\Bigg),$$
where the sets $\mathcal{O}_{\lambda,k,i}$'s are pairwise disjoint.
\end{thm}

By Theorem \ref{thmofconj1'} with Lemma \ref{lemofm-pm} and Lemma \ref{lemmaofElambdai},
we can easily give the proof of Theorem \ref{thmofconj1}.

\noindent{\emph{Proof of Theorem \ref{thmofconj1}.}
We have

\begin{align*}
  \left|\mathcal{O}_{1,k}(n)\right|&=\Bigg|\bigcup_{\lambda\in\mathcal{O}_k(n)}\Bigg(\bigcup_{i}\mathcal{O}_{\lambda,k,i}\Bigg)\Bigg|
  =\sum_{\lambda\in\mathcal{O}_k(n)}\Bigg(\sum_{i}\left|\mathcal{O}_{\lambda,k,i}\right|\Bigg)\\[6pt]
  &=\sum_{\lambda\in\mathcal{O}_k(n)}\sum_{i}\frac{m_i-p(m_i)}{k-1}=\frac{1}{k-1}\cdot\Bigg(\sum_{\lambda\in\mathcal{O}_k(n)}\ell(\lambda)-\ell(\varphi(\lambda))\Bigg)\\[6pt]
  &=\frac{1}{k-1}\cdot\Bigg(\sum_{\lambda\in\mathcal{O}_k(n)}\ell(\lambda)-\sum_{\lambda\in\mathcal{D}_k(n)}\ell(\lambda)\Bigg).
\end{align*}
\qed

\begin{table}
  \centering
  \begin{tabular}{p{2.5cm}<{\centering} p{2.5cm}<{\centering} p{4.5cm}<{\raggedright}}
  \specialrule{0.12em}{0pt}{6pt}
    $\lambda\in\mathcal{O}_3(9)$ & $\varphi(\lambda)\in\mathcal{D}_3(9)$ & \tc{$\mathcal{O}_{\lambda,3,i}\subseteq\mathcal{O}_{1,3}(9)$} \\
    \specialrule{0.12em}{3pt}{6pt}
    {$(1^{9})$} & {$(9)$} & $\{(1^6,3),(1^3,3^2),(3^3),(9)\}$\\
    \specialrule{0.05em}{3pt}{3pt}
    {$(1^7,2)$} & {$(1,2,3^2)$} & $\{(1^4,2,3),(1,2,3^2)\}$\\
    \specialrule{0.05em}{3pt}{3pt}
    {$(1^5,2^2)$} & {$(1^2,2^2,3)$} & $\{(1^2,2^2,3)\}$\\
    \specialrule{0.05em}{3pt}{3pt}
    $(1^3,2^3)$ & $(3,6)$ & $\{(2^3,3)\}$, $\{(1^3,6)\}$\\
    \specialrule{0.05em}{3pt}{3pt}
    $(1,2^4)$ & $(1,2,6)$ & $\{(1,2,6)\}$\\
    \specialrule{0.05em}{3pt}{3pt}
    $(1^5,4)$ & $(1^2,3,4)$ & $\{(1^2,3,4)\}$\\
    \specialrule{0.05em}{3pt}{3pt}
    $(1^3,2,4)$ & $(2,3,4)$ & $\{(2,3,4)\}$\\
    \specialrule{0.05em}{3pt}{3pt}
    $(1^4,5)$ & $(1,3,5)$ & $\{(1,3,5)\}$\\
    \specialrule{0.12em}{6pt}{10pt}
  \end{tabular}
  \caption{Correspondence
among the sets $\mathcal{O}_3(9)$,  $\mathcal{D}_3(9)$ and the subsets $\mathcal{O}_{\lambda,3,i}\subseteq\mathcal{O}_{1,3}(9)$}
  \label{tableofconj1}
\end{table}

\noindent{\emph{Proof of Theorem \ref{thmofconj1'}}.}
By the definition of $\mathcal{O}_{\lambda,k,i}$,
it is obvious that $\bigcup_{\lambda\in\mathcal{O}_k(n)}(\bigcup_{i}\mathcal{O}_{\lambda,k,i})\subseteq\mathcal{O}_{1,k}(n)$.
Hence to prove the theorem,
we only  need to show that for each $\pi\in\mathcal{O}_{1,k}(n)$,
there exists a unique pair $(\lambda,i)$
satisfying $\pi\in\mathcal{O}_{\lambda,k,i}$,
where $\lambda\in\mathcal{O}_k(n)$
contains $m_i$ parts $i$ with $m_i\geq k$.
Fix a partition $\pi\in\mathcal{O}_{1,k}(n)$ and suppose that the only part divisible by $k$ in $\pi$ is $a$,
we write $a$ as $a=k^j \cdot i $, where $k\nmid i$ and $j\geq1$.
Let $\lambda$ be the partition constructed from $\pi$ by replacing every part $a$  by $k^j$ parts $i$,
thus we obtain $\lambda\in\mathcal{O}_k(n)$.
Since $k^j\geq k$, $\lambda$ contains at least $k$ parts $i$,
implying $m_i\geq k$.
Therefore we find the unique required pair
$(\lambda,i)$
such that $\pi\in\mathcal{O}_{\lambda,k,i}$,
which completes the proof.
\qed

\begin{example}
  In Table \ref{tableofconj1}, we give the correspondence among the sets $\mathcal{O}_3(9)$,  $\mathcal{D}_3(9)$ and the subsets $\mathcal{O}_{\lambda,3,i}\subseteq\mathcal{O}_{1,3}(9)$,
  where each row contains a  partition $\lambda\in\mathcal{O}_3(9)$,
  the corresponding partition $\varphi(\lambda)\in\mathcal{D}_3(9)$,
  and the corresponding subsets $\mathcal{O}_{\lambda,3,i}$'s arranged in the increasing order of $i$.
  Note that for  brevity, we just list the pairs $(\lambda,\varphi(\lambda))$
  whose differences of  number of parts are nonzero.
\end{example}

\begin{remark}
  Let $\mathcal{D}_{1,k}$ denote the set of partitions of $n$ with
  exactly one part occurring at least $k$ times.
  Let $E_k(n)$ denote the difference between
  the number of parts congruent to $1$ modulo $k$ in partitions of $\mathcal{O}_k(n)$
  and the number of distinct parts of partitions in $\mathcal{D}_k(n)$.
  Fu and Tang \cite[Theorem 1.5]{Fu} analytically obtained
  that for $n\geq0$ and $k\geq2$, $\mathcal{O}_{1,k}(n)=\mathcal{D}_{1,k}(n)=E_k(n)$,
  whose combinatorial proof still remains open.
  By letting $k=2$, this theorem implies the validity of Conjecture \ref{conj1},
  where the part of $\mathcal{O}_{1,k}(n)=\mathcal{D}_{1,k}(n)$ is already ensured by Franklin's theorem.
\end{remark}

\section{Combinatorial proof of Theorem \ref{thmofconj2}}\label{secconj2}

In this section, we  introduce some new notations. 
For any positive integer $n$, denote by $f_k(n)$  the largest factor of $n$ which is not divisible by $k$.
For any positive integer $d$ with $k\nmid d$, define  $F_{k,d}=\{n\geq1\colon\,f_k(n)=d\}$.
Therefore, by the construction of Glaisher's bijection $\phi\colon\,\mathcal{D}_k(n)\rightarrow\mathcal{O}_k(n)$,
we easily obtain the following lemma.

\begin{lem}\label{lemofdifference}
For any $\lambda\in\mathcal{D}_k(n)$, we have
  \begin{equation}\label{eqoflemmaF_k}
   \bar{\ell}(\lambda)-\bar{\ell}(\phi(\lambda))=\sum_{
   k\nmid d \atop F_{k,d}\cap\lambda\neq\emptyset}\left(\left|F_{k,d}\cap \lambda\right|-1\right),
  \end{equation}
  here we view $\lambda$ as a multiset.
\end{lem}
\pf
Let $\lambda\in\mathcal{D}_k(n)$, and let $d>0$ be an integer with $k\nmid d$.
Then we know that $\bar{\ell}(\lambda)=\sum_{d}|F_{k,d}\cap\lambda|$ and  $\bar{\ell}(\phi(\lambda))=\sum_{d}|F_{k,d}\cap\phi(\lambda)|$
since $\{F_{k,d}\}$ is a
partition of the positive integer set.
By the definition of Glaisher's bijection $\phi$, $|F_{k,d}\cap\phi(\lambda)|\neq\emptyset$
if and only if $|F_{k,d}\cap\lambda|\neq\emptyset$.
Thus, the difference between $\bar{\ell}(\lambda)$ and $\bar{\ell}(\phi(\lambda))$
is the sum of the differences between $|F_{k,d}\cap\lambda|$
and $|F_{k,d}\cap\phi(\lambda)|$
for any $|F_{k,d}\cap\lambda|\neq\emptyset$.
Therefore \eqref{eqoflemmaF_k} holds
since $F_{k,d}\cap\phi(\lambda)\in\{\{d\},\emptyset\}$
by Glaisher's bijection.
\qed


Let $\lambda\in\mathcal{D}_k(n)$, for any $F_{k,d}\cap \lambda\neq \emptyset$ with $d$, $k\nmid d$,
we construct $\mathcal{T}_{\lambda,k,d}\subseteq\mathcal{T}_k(n)$
of size $|F_{k,d}\cap\lambda|-1$.
Assuming $|F_{k,d}\cap\lambda|=p$, then $\lambda$ contains
\begin{equation*}
  \big\{\,\underbrace{k^{a_1}d,k^{a_1}d,\ldots,k^{a_1}d}_{\leq k-1},\,\underbrace{k^{a_2}d,k^{a_2}d,\ldots,k^{a_2}d}_{\leq k-1},\ldots,\underbrace{k^{a_p}d,k^{a_p}d,\ldots,k^{a_p}d}_{\leq k-1}\,\big\}
\end{equation*}
as parts with $0\leq a_1<a_2<\cdots<a_p$.
For $1\leq i\leq p-1$,  we construct $\tau_d^i\in\mathcal{T}_{\lambda,k,d}$ from $\lambda$ by replacing
one part  $k^{a_{i+1}}d$ in $\lambda$ by $k$ parts  $k^{a_i}d$,
$k-1$ parts  $k^{a_i+1}d$, $k-1$ parts $k^{a_i+2}d$, $\ldots$ ,
$k-1$ parts  $k^{a_{i+1}-1}d$.
Since
\[
k\cdot k^{a_i}d+(k-1)\cdot k^{a_i+1}d+(k-1)\cdot k^{a_i+2}d+\cdots +(k-1)\cdot k^{a_{i+1}-1}d=k^{a_{i+1}}d,
\]
it follows that $|\tau^i_d|=|\lambda|$.
Moreover, the number of parts $k^{a_i}d$ in $\tau_d^i$ is between $k+1$ and $2k-1$
while other parts in $\tau^i_d$ occur less than $k$ times,
therefore $\tau_d^i\in\mathcal{T}_k(n)$.
Define $$\mathcal{T}_{\lambda,k,d}=\left\{\tau^1_d,\tau^2_d,\ldots,\tau^{p-1}_d\right\}.$$
Note that $\mathcal{T}_{\lambda,k,d}=\emptyset$ when $|F_{k,d}\cap\lambda|=1$.

\begin{thm}\label{thmofconj2'}
For $k\geq1$ and $n\geq0$, we have
  \[\mathcal{T}_k(n)=\bigcup_{\lambda\in\mathcal{D}_k(n)}\Bigg(\bigcup_{k\nmid d \atop F_{k,d}\cap\lambda\neq\emptyset}\mathcal{T}_{\lambda,k,d}\Bigg),\]
where the sets $\mathcal{T}_{\lambda,k,d}$'s are pairwise disjoint.
\end{thm}

Therefore by the combination of Theorem \ref{thmofconj2'} and Lemma \ref{lemofdifference}, we can give the proof
of Theorem \ref{thmofconj2}.

\noindent{\emph{Proof of Theorem \ref{thmofconj2}.}}
We obtain that
\begin{align*}
  \left|\mathcal{T}_k(n)\right|&=\sum_{\lambda\in\mathcal{D}_k(n)}
  \sum_{k\nmid d \atop F_{k,d}\cap\lambda\neq\emptyset}\left|\mathcal{T}_{\lambda,k,d}(n)\right|=\sum_{\lambda\in\mathcal{D}_k(n)}
  \sum_{k\nmid d \atop F_{k,d}\cap\lambda\neq\emptyset}\left(\left|F_{k,d}\cap\lambda\right|-1\right)\\[6pt]
  &=\sum_{\lambda\in\mathcal{D}_k(n)}\big(\bar{\ell}(\lambda)-\bar{\ell}(\phi(\lambda))\big)
  =\sum_{\lambda\in\mathcal{D}_k(n)}\bar{\ell}(\lambda)-\sum_{\lambda\in\mathcal{O}_k(n)}\bar{\ell}(\lambda).
\end{align*}}
\qed

\begin{table}
  \centering
  \begin{tabular}{p{2.5cm}<{\centering} p{3cm}<{\centering} p{4cm}<{\centering}}
  \specialrule{0.12em}{0pt}{6pt}
    $\lambda\in\mathcal{D}_3(12)$ & $\phi(\lambda)\in\mathcal{O}_3(12)$ & $\mathcal{T}_{\lambda,3,d}\subseteq\mathcal{T}_3(12)$ \\
    \specialrule{0.12em}{3pt}{6pt}
    $(1^2,2^2,3^2)$ & $(1^8,2^2)$ & $\{(1^5,2^2,3)\}$ \\
    \specialrule{0.05em}{3pt}{3pt}
    $(1,2^2,3,4)$  & $(1^4,2^2,4)$ & $\{(1^4,2^2,4)\}$ \\
    \specialrule{0.05em}{3pt}{3pt}
    $(1^2,3^2,4)$ & $(1^8,4)$ & $\{(1^5,3,4)\}$\\
    \specialrule{0.05em}{3pt}{3pt}
    $(1,3,4^2)$ & $(1^4,4^2)$ & $\{(1^4,4^2)\}$\\
    \specialrule{0.05em}{3pt}{3pt}
    $(1^2,2,3,5)$ & $(1^5,2,5)$ & $\{(1^5,2,5)\}$\\
    \specialrule{0.05em}{3pt}{3pt}
    $(1,3^2,5)$ & $(1^7,5)$ & $\{(1^4,3,5)\}$\\
    \specialrule{0.05em}{3pt}{3pt}
    $(1^2,2^2,6)$ & $(1^2,2^5)$ & $\{(1^2,2^5)\}$\\
    \specialrule{0.05em}{3pt}{3pt}
    $(1,2,3,6)$ & $(1^4,2^4)$ & $\{(1^4,2,6)\}$, $\{(1,2^4,3)\}$\\
    \specialrule{0.05em}{3pt}{3pt}
    $(2,4,6)$ & $(2^4,4)$ & $\{(2^4,4)\}$\\
    \specialrule{0.05em}{3pt}{3pt}
    $(1^2,3,7)$ & $(1^5,7)$ & $\{(1^5,7)\}$\\
    \specialrule{0.05em}{3pt}{3pt}
    $(1,3,8)$ & $(1^4,8)$ & $\{(1^4,8)\}$\\
    \specialrule{0.05em}{3pt}{3pt}
    $(1,2,9)$ & $(1^{10},2)$ & $\{(1^{4},2,3^2)\}$\\
    \specialrule{0.05em}{3pt}{3pt}
    $(3,9)$ & $(1^{12})$& $\{(3^4)\}$\\
    \specialrule{0.12em}{6pt}{10pt}
  \end{tabular}
  \caption{Correspondence
among the sets $\mathcal{D}_3(12)$, $\mathcal{O}_3(12)$  and the subsets $\mathcal{T}_{\lambda,3,d}\subseteq\mathcal{T}_3(12)$\label{tableofconj2}}
\end{table}

\noindent{\emph{Proof of Theorem \ref{thmofconj2'}}.}
By the construction of $\mathcal{T}_{\lambda,k,d}$,
 it is clear that if we fix any $\lambda\in\mathcal{D}_k(n)$,
then $\mathcal{T}_{\lambda,k,d_1}\cap\mathcal{T}_{\lambda,k,d_2}=\emptyset$
if $d_1\neq d_2$
since $F_{k,d_1}\cap F_{k,d_2}=\emptyset$.
Hence to complete the proof,
we only need to show that
for arbitrary partition $\tau\in\mathcal{T}_k(n)$,
there  exists one and only one partition $\lambda\in\mathcal{D}_k(n)$
such that $\tau\in\mathcal{T}_{\lambda,k,d}$
for some  $d$, $k\nmid d$.

Assume $\tau\in\mathcal{T}_k(n)$
and there is only one part $k^ad$ in $\lambda$ for some $a\geq 0$,
the multiplicity of which is more than $k$ and less than $2k$.
We replace $k$ parts $k^ad$ by one part $k^{a+1}d$,
which preserves the weight of partition and
reduces the  multiplicity of part $k^ad$
to at most $k-1$.
If it causes the multiplicity of part $k^{a+1}d$  increased to $k$,
which is possible,
then we continue to replace $k$ parts $k^{a+1}d$ by one part $k^{a+2}d$.
Thus as long as the multiplicity of part $k^{a+i}d$ is $k$,
we replace these $k$ parts $k^{a+i}d$ by one part $k^{a+i+1}d$.
This process must stop at $a+m$ for some $m\geq1$
since the number of part of $\lambda$ is finite.
Then we obtain the desired $\lambda\in\mathcal{D}_k(n)$ since
the number of part $k^{a+i}d$ is less than $k$ for all $0\leq i\leq m$
and other parts of $\tau$  whose
multiplicities are less than $k$ originally
stay unchanged.
Note that $\{k^ad,k^{a+1}d,\ldots,k^{a+m}d\}\subseteq F_{k,d}\cap\lambda$, which implies $F_{k,d}\cap\lambda\neq\emptyset$.
Therefore $\tau\in\mathcal{T}_{\lambda,k,d}$ and the proof is completed.
\qed

\begin{example}
  We give the correspondence among the sets $\mathcal{D}_3(12)$, $\mathcal{O}_3(12)$  and the subsets $\mathcal{T}_{\lambda,3,d}\subseteq\mathcal{T}_3(12)$ in Table \ref{tableofconj2},
  where each row contains a partition $\lambda\in\mathcal{D}_3(12)$,
  the corresponding  $\phi(\lambda)\in\mathcal{O}_3(12)$,
  and the  corresponding subsets $\mathcal{T}_{\lambda,3,d}$'s arranged in the increasing order of $d$.
  For the sake of conciseness,
  we just list the pairs $(\lambda,\phi(\lambda))$
  whose differences of the number of distinct parts are nonzero.
\end{example}

\section{Combinatorial proof of Conjecture \ref{conj3}}\label{secconj3}

In this section, we  prove Conjecture \ref{conj3} which is restated below.
\begin{conjprime}{conj3}\label{conj3'}
Let $\mathcal{D_O}(n)$ be the set of distinct partitions of  $n$
with odd length,
and $\mathcal{G}(n)$ be the set of
gap-free partitions of $n$.
Then
\begin{equation*}
|\mathcal{G}(n)|=\sum_{\lambda\in\mathcal{D_O}(n)}s(\lambda),
\end{equation*}
where $s(\lambda)$ is the smallest part of $\lambda$.
\end{conjprime}


Denote by $\mathcal{G_I}(n)$ the set of  gap-free partitions of $n$
with the smallest part  1 and the largest part  odd.
For any partition $\lambda$,
let  $m(\lambda)$ be the multiplicity of the largest part of $\lambda$.
Then we have the following lemma.
\begin{lem}\label{D(n)toG(n)}
  The conjugate of a partition gives a bijection between $\mathcal{G}_\mathcal{I}(n)$ and $\mathcal{D}_\mathcal{O}(n)$.
  Thus for any $\lambda\in\mathcal{G}_\mathcal{I}(n)$, we have
  $\lambda'\in\mathcal{D}_\mathcal{O}(n)$ and
  $m(\lambda)=s(\lambda')$.
\end{lem}

%
%
%
%
%

For convenience, we need to introduce
some notations and operators on partitions.
Let $\lambda=(1^{m_1},2^{m_2},\ldots,k^{m_k})$
be a partition
with $m_i\geq 0$ for $1\leq i\leq k$.
We call every $i^{m_i}$ the \emph{block} of $i$ if $m_i\geq 1$
and denoted by $B_i$.
Define $B_i<B_j$ if and only if $i<j$,
and $B_i$ to be  odd (resp. even) if and only if $i$ is odd (resp. even).
The number of blocks of $\lambda$ is exactly the number of distinct parts in $\lambda$.
Hence if a partition $\lambda$ is gap-free,
then $\lambda$ contains  $\bar{\ell}(\lambda)$   consecutively-indexed blocks,
and $\mathcal{G_I}$  consists of the gap-free partitions
with the smallest block  $B_1$ and the largest block  odd.
Let $\lambda$ be a gap-free partition of the form
\begin{equation*}
  \Big(\,\underbrace{k,\ldots,\bm{k}}_{B_k},\,\underbrace{k+1,\ldots, \bm{k+1}}_{B_{k+1}},\ldots,\underbrace{k+r-1,\ldots, \bm{k+r-1}}_{B_{k+r-1}},\ldots,\underbrace{k+\ell-1,\ldots, \bm{k+\ell-1}}_{B_{k+\ell-1}}\,\Big),
\end{equation*}
then for any $0\leq r\leq\ell$, we can define an \emph{increasing} operator $\xi_r$
which can increase the weight of $\lambda$ by $r$ and
preserve its property of being gap-free.
To be more specific, notice that $B_k,B_{k+1},\ldots,B_{k+r-1}$ are the $r$ smallest blocks of $\lambda$,
then for each $1\leq j\leq r$,
add 1 to the last part $k+j-1$ in  block $B_{k+j-1}$,
which is in bold type.
Let $\xi_r(\lambda)$ be the resulting partition, then
it is also gap-free and of weight $|\lambda|+r$.
Symmetrically we can define  a \emph{decreasing} operator $\xi^-_r$, where $0\leq r\leq \ell$.
Note that  $B_{\ell-r+1},B_{\ell-r+2},\ldots,B_{\ell}$
are the $r$ largest blocks of $\lambda$.
$\xi^{-}_r(\lambda)$ is obtained from $\lambda$
by subtracting 1 from the first part $\ell-j+1$ in  block $B_{\ell-j+1}$
for each $1\leq j\leq r$.
Thus $\xi^{-}_r(\lambda)$ is a gap-free partition of weight $|\lambda|-r$.

\begin{lem}\label{lemmaofxi}
  Let $\lambda$ be
  a gap-free partition with $s(\lambda)>1$ and $r\in\{\bar{\ell}(\lambda)-1,\bar{\ell}(\lambda)\}$,
  then the following hold:
  \begin{enumerate}[(1).]
    \item \label{lemmaofxi1}
    $\bar{\ell}(\xi_r(\lambda))=r$ or $\bar{\ell}(\xi_r(\lambda))=r+1$;
    \item \label{lemmaofxi2}
    $\bar{\ell}(\xi^-_r(\lambda))=r$ or $\bar{\ell}(\xi^-_r(\lambda))=r+1$;
    \item \label{lemmaofxi3}
    $\xi_r^{-1}=\xi^-_r$, where $\xi_r^{-1}$ is the inverse operator of $\xi_r$.
  \end{enumerate}
\end{lem}
\pf
Let $\lambda=(k^{m_k},(k+1)^{m_{k+1}},\ldots,(k+\ell-1)^{m_{k+\ell-1}})$ be the gap-free partition
with $k>1$ and $\bar{\ell}(\lambda)=\ell$.
Set $r=\ell$, then $\xi_{\ell}(\lambda)=\mu=(k^{m_k-1},(k+1)^{m_{k+1}},\ldots,(k+\ell-1)^{m_{k+\ell-1}},k+\ell)$.
Note that $m_k-1$ can be zero, which means that $\bar{\ell}(\mu)=\ell$ or $\bar{\ell}(\mu)=\ell+1$.
Hence the action of $\xi^-_\ell$ is valid
and $\xi^-_\ell(\mu)=(k^{m_k},(k+1)^{m_{k+1}},\ldots,(k+\ell-1)^{m_{k+\ell-1}})=\lambda$.
On the other hand,
$\xi^-_\ell(\lambda)=\nu=(k-1,k^{m_k},\ldots,(k+\ell-1)^{m_{k+\ell-1}-1})$,
implying that $\bar{\ell}(\nu)=\ell$ or $\bar{\ell}(\nu)=\ell+1$.
It follows that $\xi_\ell(\nu)=(k^{m_k},(k+1)^{m_{k+1}},\ldots,(k+\ell-1)^{m_{k+\ell-1}})=\lambda$.
Thus we have verified  the lemma for the case  $r=\bar{\ell}(\lambda)$.
The proof for the case  $r=\bar{\ell}(\lambda)-1$
is completely the same so we leave it to the readers.
\qed

Let $\lambda=(1^{m_1},2^{m_2},\ldots,\ell^{m_\ell})\in\mathcal{G}(n)$
be a gap-free partition with the smallest block $B_1$ and $\bar{\ell}(\lambda)=\ell$.
Note that here we do not require the constraint that $B_\ell$ is odd.
Then we shall define a series of operators $\varrho_i$ for $1\leq i\leq m_\ell$ as follows.
First we delete $i-1$ parts $\ell$ from $\lambda$ then $\lambda$ becomes
$\lambda^{i}_1=(1^{m_1},2^{m_2},\ldots,\ell^{m_\ell-i+1})$.
It is clear that $\lambda^{i}_1$ is also gap-free and  $\bar{\ell}(\lambda^{i}_1)=\ell$.
For $2\leq j\leq i$, let $\lambda^i_j=\xi_\ell(\lambda^i_{j-1})$
iteratively and finally let $\varrho_i(\lambda)=\lambda^i_{i}$.
It is easy to see that $\varrho_i(\lambda)$ and $\varrho_j(\lambda)$
are different whenever $i\neq j$.


\begin{thm}\label{thmofconj3}
  For any $\lambda\in\mathcal{G_I}(n)$,
  define $\mathcal{G}_\lambda=\{\varrho_i(\lambda):1\leq i\leq m(\lambda)\}.$
  Then the sets $\mathcal{G}_\lambda$'s are pairwise disjoint and
  $$\mathcal{G}(n)=\bigcup_{\lambda\in\mathcal{G_I}(n)}\mathcal{G}_\lambda.$$
\end{thm}

Combining  Theorem \ref{thmofconj3} and Lemma \ref{D(n)toG(n)},
we see that $$|\mathcal{G}(n)|=\bigg|\bigcup_{\lambda\in\mathcal{G_I}(n)}\mathcal{G}_\lambda\bigg|
=\sum_{\lambda\in\mathcal{G_I}(n)}|\mathcal{G}_\lambda|=\sum_{\lambda\in\mathcal{G_I}(n)}m(\lambda)
=\sum_{\lambda\in\mathcal{D_O}(n)}s(\lambda),$$
which confirms Conjecture \ref{conj3'}.

\noindent{\emph{Proof of Theorem \ref{thmofconj3}.}}
We first clarify that the operators $\varrho_i$'s are well-defined.
To this end, by the definition of increasing operator $\xi_r$ for $r\geq 0$,
we only need to show that
when we apply $\xi_r$ on some partition $\lambda$, the number of the blocks of $\lambda$ is no less than $r$.
Actually, we can claim that
for any $\lambda\in\mathcal{G_I}(n)$ and $1\leq i \leq m(\lambda)$,
we have
$\bar{\ell}(\lambda^i_j)=\bar{\ell}(\lambda)$
or $\bar{\ell}(\lambda^i_j)=\bar{\ell}(\lambda)+1$
for $1\leq j\leq i$.
Let $\lambda\in\mathcal{G_I}(n)$ be a gap-free partition with $\bar{\ell}(\lambda)=\ell$ and $m(\lambda)=m$,
then for any $1\leq i\leq m$, it is evident that  $\bar{\ell}(\lambda^i_1)=\ell$.
By  induction, we assume that $\bar{\ell}(\lambda^i_j)=\ell$ or
$\bar{\ell}(\lambda^i_j)=\ell+1$ for  $2\leq j\leq i-1$.
By applying  $\xi_\ell$ on $\lambda^i_j$, it produces $\xi_\ell(\lambda^i_j)=\lambda^i_{j+1}$.
Since $\ell\in\{\bar{\ell}(\lambda^i_j)-1,\bar{\ell}(\lambda^i_j)\}$\
by the induction hypothesis,
then by Lemma \ref{lemmaofxi} \eqref{lemmaofxi1},
we have $\bar{\ell}(\lambda^i_{j+1})=\ell$ or
$\bar{\ell}(\lambda^i_{j+1})=\ell+1$.
Hence we derive that $\bar{\ell}(\lambda^i_j)=\bar{\ell}(\lambda)$
or $\bar{\ell}(\lambda^i_j)=\bar{\ell}(\lambda)+1$
for $1\leq j\leq i$.
 Thus the operators   $\varrho^i$'s are well-defined.

Next, for every  $\lambda\in\mathcal{G_I}(n)$, we shall show that $\mathcal{G}_\lambda\subseteq\mathcal{G}(n)$, which implies that  $\bigcup_{\lambda\in\mathcal{G_I}(n)}\mathcal{G}_\lambda\subseteq\mathcal{G}(n)$.
Since for any $\lambda\in\mathcal{G_I}(n)$
and $1\leq i\leq m(\lambda)$,
$\lambda^i_1$ is gap-free of weight $n-(i-1)\bar{\ell}(\lambda)$
and each time applying $\xi_{\bar{\ell}(\lambda)}$ preserves the property of being gap-free
and increases the weight  by $\bar{\ell}(\lambda)$,
then after $(i-1)$-time composition of $\xi_{\bar{\ell}(\lambda)}$,
it is clear that
$\varrho_i(\lambda)$ is gap-free and of weight $n$.
Thus we have $\mathcal{G}_\lambda\subseteq\mathcal{G}(n)$
then $\bigcup_{\lambda\in\mathcal{G_I}(n)}\mathcal{G}_\lambda\subseteq\mathcal{G}(n)$.

Finally, we prove that for any $\mu\in \mathcal{G}(n)$,
there  exists only one $\lambda\in\mathcal{G_I}(n)$
such that $\mu\in\mathcal{G}_\lambda$.
Given any $\mu\in\mathcal{G}(n)$, if $\mu\in\mathcal{G_I}(n)$, then we are done;
otherwise let $\mu=(k^{m_k},(k+1)^{m_{k+1}},\ldots,(k+\ell-1)^{m_{k+\ell-1}})$
with $\bar{\ell}(\mu)=\ell$.
Set $\mu_1=\mu$ and $r=2\lceil\ell/2\rceil-1$, i.e. $r$ is the largest odd number that do not exceed  $\ell$,
so that $\bar{\ell}(\mu_1)=r$ or $\bar{\ell}(\mu_1)=r+1$.
For $i\geq 2$,
we keep constructing $\mu_i$ by setting $\mu_i=\xi^-_r(\mu_ {i-1})$
until we first get $\mu_{m}\in\mathcal{G_I}$ for some $m\geq 2$.
This can be done since for $\ell\geq 1$, we have $r\geq 1$ and for $1\leq i\leq m$, $\bar{\ell}(\mu_i)=r$ or $\bar{\ell}(\mu_i)=r+1$ by Lemma \ref{lemmaofxi} \eqref{lemmaofxi2}.
Since $\xi^-_r$ is the inverse operator of $\xi_r$ by Lemma \ref{lemmaofxi} \eqref{lemmaofxi3},
it can be checked that the largest block of $\mu_m$ is $B_r$.
Hence, by adding $m-1$ parts  $r$ back to $\mu_m$, we obtain the unique $\lambda\in\mathcal{G_I}(n)$.
This completes the proof.
\qed

\begin{example}
  We list  the sets $\mathcal{D_O}(12)$, $\mathcal{G_I}(12)$ and the subsets $\mathcal{G}_\lambda\subseteq\mathcal{G}(12)$ in Table \ref{tableofconj3},
  where each row contains a gap-free partition $\lambda\in\mathcal{G_I}(12)$,
  the corresponding distinct partition
 $\lambda'\in\mathcal{D_O}(12)$,
  and the corresponding subset $\mathcal{G}_\lambda$ of which
  partitions  are arranged from $\varrho_1(\lambda)$ to $\varrho_{m(\lambda)}(\lambda)$.
\end{example}
\begin{table}[h!]
  \centering
  \begin{tabular}{p{3cm}<{\centering} p{3cm}<{\centering} p{5.5cm}<{\raggedright}}
  \specialrule{0.12em}{0pt}{6pt}
    $\lambda'\in\mathcal{D_O}(12)$ & $\lambda\in\mathcal{G_I}(12)$ & \tc{$\mathcal{G}_\lambda\subseteq\mathcal{G}(12)$} \\
    \specialrule{0.12em}{3pt}{6pt}
    {(5,4,3)} & {$(1,2,3^3)$} & $\{(1,2,3^3),(2,3^2,4),(3,4,5)\}$ \\
    \specialrule{0.05em}{3pt}{3pt}
    {(6,4,2)} & {$(1^2,2^2,3^2)$} & $\{(1^2,2^2,3^2),(1,2^2,3,4)\}$\\
    \specialrule{0.05em}{3pt}{3pt}
    {(6,5,1)} & {$(1,2^4,3)$} & $\{(1,2^4,3)\}$\\
    \specialrule{0.05em}{3pt}{3pt}
    (7,3,2) & $(1^4,2,3^2)$ & $\{(1^4,2,3^2),(1^3,2,3,4)\}$\\
    \specialrule{0.05em}{3pt}{3pt}
    (7,4,1) & $(1^3,2^3,3)$ & $\{(1^3,2^3,3)\}$\\
    \specialrule{0.05em}{3pt}{3pt}
    (8,3,1) & $(1^5,2^2,3)$ & $\{(1^5,2^2,3)\}$\\
    \specialrule{0.05em}{3pt}{3pt}
    (9,2,1) & $(1^7,2,3)$ & $\{(1^7,2,3)\}$\\
    \specialrule{0.05em}{3pt}{3pt}
    \multirow{3}{*}{(12)} & \multirow{3}{*}{$(1^{12})$} & $\{(1^{12}),(1^{10},2),(1^8,2^2),(1^6,2^3),$\\[1.5pt]
    & & $(1^4,2^4),(1^2,2^5),(2^6),(2^3,3^2),$\\[1.5pt]
    & & $(3^4),(4^3),(6^2),(12)\}$\\
    \specialrule{0.12em}{6pt}{10pt}
  \end{tabular}

  \caption{Correspondence among the sets $\mathcal{D_O}(12)$, $\mathcal{G_I}(12)$ and the subsets $\mathcal{G}_\lambda\subseteq\mathcal{G}(12)$}
  \label{tableofconj3}
\end{table}

Let $\mathcal{D_E}(n)$ be the set of distinct partitions of $n$ with even length.
We can also investigate the quantitative relationship between
the sets $\mathcal{D_E}(n)$ and  $\mathcal{G}(n)$.
Denote by $\mathcal{G_I'}(n)$ the set of
gap-free partitions of $n$
with the smallest block  $B_1$
and the largest block even.
By Lemma \ref{D(n)toG(n)}, we know that
taking conjugate also gives a bijection
from $\mathcal{G_I'}(n)$ to $\mathcal{D_E}(n)$ such that
$m(\lambda)=s(\lambda')$ for any $\lambda\in\mathcal{G_I'}(n)$.

\begin{thm}\label{thmofconj3'}
For any $\lambda\in\mathcal{G_I'}(n)$,
define $\mathcal{G}_\lambda=\{\varrho_i(\lambda):1\leq i\leq m(\lambda)\}$ as before.
Then the sets $\mathcal{G}_\lambda$'s are pairwise disjoint and
  $$\mathcal{G}(n)\left\backslash\mathcal{G}^0(n)\right.=\bigcup_{\lambda\in\mathcal{G_I'}(n)}\mathcal{G}_\lambda,$$
  where $\mathcal{G}^0(n)$ is the set of gap-free partitions of $n$ with only one block.
\end{thm}

\pf The proof is exactly the same as Theorem \ref{thmofconj3}
 except for setting $r=2\lfloor\ell/2\rfloor$, i.e., $r$ is the largest even number
 that do not exceed $\ell$.
\qed

Hence we have proved the following theorem.
\begin{thm}
The number of gap-free partitions of $n$ with at least two blocks is equal to the sum of the smallest parts in distinct partitions of $n$ with an even number of parts:
$$\vert \mathcal{G}(n)\setminus \mathcal{G}^{0}(n) \vert=\sum_{\lambda\in \mathcal{D_E}(n)}s(\lambda).$$
\end{thm}

\begin{example}
  We list the sets $\mathcal{D_E}(12)$, $\mathcal{G_I'}(12)$ and the subsets $\mathcal{G}_\lambda\subseteq\mathcal{G}(12)\left\backslash\mathcal{G}^0(12)\right.$ in Table \ref{tableofconj3'},
  where each row contains a gap-free partition $\lambda\in\mathcal{G_I'}(12)$,
  the corresponding  distinct partition  $\lambda'\in\mathcal{D_E}(12)$,
  and the corresponding subset $\mathcal{G}_\lambda$ of which
  partitions  are arranged from $\varrho_1(\lambda)$ to $\varrho_{m(\lambda)}(\lambda)$.
\end{example}

\begin{table}[h!]
  \centering
  \begin{tabular}{p{3cm}<{\centering} p{3cm}<{\centering} p{5.5cm}<{\raggedright}}
  \specialrule{0.12em}{0pt}{6pt}
    $\lambda'\in\mathcal{D_E}(12)$ & $\lambda\in\mathcal{G_I'}(12)$ & \tc{$\mathcal{G}_\lambda\subseteq\mathcal{G}(12)\left\backslash\mathcal{G}^0(12)\right.$} \\
    \specialrule{0.12em}{3pt}{6pt}
    {(5,4,2,1)} & {$(1,2^2,3,4)$} & $\{(1,2^2,3,4)\}$ \\
    \specialrule{0.05em}{3pt}{3pt}
    {(6,3,2,1)} & {$(1^3,2,3,4)$} & $\{(1^3,2,3,4)\}$\\
    \specialrule{0.05em}{3pt}{3pt}
    \multirow{2}{*}{(7,5)} & \multirow{2}{*}{$(1^2,2^5)$} & \{$(1^2,2^5)$, $(1,2^4,3)$, $ (2^3,3^2)$, \\[1.5pt]
    & & $(2,3^2,4)$, $(3,4,5)$\}\\
    \specialrule{0.05em}{3pt}{3pt}
    \multirow{2}{*}{(8,4)} & \multirow{2}{*}{$(1^4,2^4)$} & \{$(1^4,2^4)$, $(1^3,2^3,3)$, $ (1^2,2^2,3^2)$, \\[1.5pt]
    & & $(1,2,3^3)$\}\\
    \specialrule{0.05em}{3pt}{3pt}
    {(9,3)} & {$(1^6,2^3)$} & \{$(1^6,2^3)$, $(1^5,2^2,3)$, $(1^4,2,3^2)$\} \\
    \specialrule{0.05em}{3pt}{3pt}
    {(10,2)} & {$(1^8,2^2)$} & \{$(1^8,2^2)$, $(1^7,2^2,3)$\}\\
    \specialrule{0.05em}{3pt}{3pt}
    {(11,1)} & {$(1^{10},2)$} & \{$(1^{10},2)$\}\\
    \specialrule{0.12em}{6pt}{10pt}
  \end{tabular}

  \caption{Correspondence among the sets $\mathcal{D_E}(12)$, $\mathcal{G_I'}(12)$ and the subsets $\mathcal{G}_\lambda\subseteq\mathcal{G}(12)\left\backslash\mathcal{G}^0(12)\right.$}
  \label{tableofconj3'}
\end{table}

\begin{remark}
  Note that the number of gap-free partitions of $n$ with only one block equals the number of divisors of $n$,
  which is usually denoted by $d(n)$ in  number theory and combinatorics. Thus by  Theorem
  \ref{thmofconj3} and Theorem \ref{thmofconj3'}, we rediscover the following identity
  \begin{equation}\label{divisorsofn}
    \sum_{\lambda\in\mathcal{D_O}(n)}s(\lambda)-\sum_{\lambda\in\mathcal{D_E}(n)}s(\lambda)=d(n),
  \end{equation}
which was first obtained by Fokkink, Fokkink and Wang  \cite{Fokkink-Fokkink-Wang}
via the sum of polynomial quotients for $n$.
In \cite{Andrews08}, Andrews asserted that \eqref{divisorsofn} is a corollary of the differentiation of the $q$-analog of Gauss's theorem \cite[ Corollary 2.4]{Andrews98} then
gave the analytic proof of \eqref{divisorsofn} by $q$-series.
\end{remark}

\noindent{\bf Acknowledgements.}
The author appreciates the anonymous referees for their helpful comments
and language editing which have greatly
improved the quality of this manuscript.
This work was supported by Scientific and Technological Research Program of Chongqing Municipal Education Commission (Grant No. KJQN201800604), 
and Doctor Scientific Research Foundation of Chongqing University of Posts and Telecommunications (Grant No. A2017-123).

\end{document}